\documentclass[a4paper,10pt]{article}

\title{ The linearity problem for the unitriangular automorphism groups of free groups\footnote{The investigation was supported by The Ministry of  Education and Science of Russian Federation, projects 14.B37.21.0359/0859.}}
\author{V. Roman'kov}

\date{}

\usepackage{amsmath,amsthm,amsfonts}
\usepackage{graphicx}
\usepackage{hyperref}
\usepackage[usenames]{color}

\newtheorem{theorem}{Theorem}[section]

\theoremstyle{definition}

\newcounter{comcount}

\begin{document}

\maketitle

\begin{abstract}
We prove that the unitriangular automorphism group of a free group of rank  $n$ has a faithful representation by matrices over a field, or in other words,  it is a linear group, if and only if  $n \leq 3.$ Thus, we have completed a description of relatively free groups  with linear the unitriangular automorphism groups. This description was 
 initiated by Erofeev and the author in \cite{Erofeev}, where  proper varieties of groups have been considered.
\end{abstract} 

\section{Introduction}

\medskip
\noindent
For each positive integer  $n,$ let   $F_n$ be a free group of rank $n$  with basis (in other words, free generating set)   $\{f_1, ... , f_n\}.$   For any     $m \leq n$   $F_m$ is considered as   subgroup gp$(f_1, ... , f_m)$ of $F_n.$ 
For any variety of groups ${\cal G},$ let  ${\cal G}(F_n)$ denote the verbal subgroup of $F_n$ corresponding to  ${\cal G}.$  Let   $G_n = F_n/{\cal G}(F_n).$  Then $G_n$ is a relatively free group of rank   $n$ in the variety  ${\cal G}.$
By  basis of $G_n$ we mean a subset $S$ such that every map of $S$ into $G_n$ extends, uniquely, to an endomorphism of $G_n.$ Write $\bar{f}_i = f_i{\cal G}(F_n)$ for $i = 1, ... , n.$ Then   $\bar{f}_1, ... , \bar{f}_n$ is a basis of $G_n.$  

Let  ${\cal G}$ be a variety of groups. Let                $G_n$ be the relatively free group corresponding to  ${\cal G}$ with basis  $\{f_1, ... , f_n\}.$  For any
    $m \leq n$    $G_m$ is considered as   subgroup gp$(f_1, ... , f_m)$ of       $G_n.$ An automorphism $\varphi  $  of      $G_n$   is called {\it unitriangular} (w.r.t. the given basis)   if  $\varphi $ is  defined by a map of the form:

\begin{equation}
\label{eq:1}
\varphi : f_1 \mapsto f_1, f_i \mapsto u_if_i \   {\textrm{для}} \ i = 2, ... , n,
\end{equation}

\noindent
where  $u_i = u_i(f_1, ... , f_{i-1})$ is an element of  $G_{i-1}.$ Every tuple of elements  $(u_2, ... , u_n)$ with this condition defines, uniquely,  automorphism of  $G_n.$ 
Let  $U_n$ be  subgroup consisting of all unitriangular (w.r.t. a given basis) automorphisms of  $G_n.$ Then it is called {\it the unitriangular automorphism group} of  $G_n.$    As  abstract group $U_n$ does not depend of a basis.   

The question of linearity of $U_n$ for an arbitrary proper  variety ${\cal G}$ has been studied by Erofeev and the author in   \cite{Erofeev}.  All cases of linearity of $G_n$ have been described.  The following Section 1 contains this description. Also, we observe some relative results on linearity for relatively free groups and algebras. 

In this paper we study the only open after \cite{Erofeev} case when   ${\cal G}$ is the variety of all groups. 
Our main result is given by the following theorem. 

\begin{theorem}
\label{th:1.1}
 \textit{The group $U_n$ of unitriangular automorphisms  of the free group $F_n$  of rank  $n$ is linear if and only if  $n \leq 3.$}
\end{theorem}

\medskip
Hence, we  complete a description of all  cases when the unitriangular automorphism group  $U_n,$ corresponding to an arbitrary variety of groups  ${\cal G},$ including the variety of all groups, is linear.

\section{Some results on linearity.}

\medskip
\noindent
We observe some results concerning the linearity of the automorphism groups and their subgroups of  relatively free groups and algebras. Recall that group $G$ is said to be 
{\it virtually nilpotent} if it has a nilpotent subgroup of  finite index. 

The linearity of Aut$(F_2)$ follows by \cite{Dyer}     from the linearity 
of the $4-$string braid group $B_4$, which is due to Krammer   \cite{Krammer}.  Bigelow \cite{Bigelow} and also Krammer \cite{Krammer1} determined that the braid group $B_n$ is linear for every $n.$   Formanek and Procesi  in \cite{Formanek} have demonstrated that Aut$(F_n)$ is not linear for $n \geq 3.$ 

Auslander and Baumslag \cite{Auslender} determined that for every finitely generated virtually nilpotent group $G$ the automorphism group   Aut$(G)$ is linear. 
Moreover, Aut$(G)$ has a faithful matrix representation over the integers 
 ${\mathbb{Z}}.$  In particular, for every relatively free virtually nilpotent group $G_n,$ the automorphism group Aut$(G_n)$  is linear over ${\mathbb{Z}}.$

\medskip
Olshanskii \cite{Olshanskii} proved for any  relatively free group  $G_n,$ which is not virtually nilpotent and is not free,  that the automorphism group  Aut$(G_n)$ is not linear. His approach does not give an information on the linearity of the unitriangular automorphism groups  $U_n$ for such relatively free groups  $G_n.$ 

Erofeev and the author \cite{Erofeev} proved for every proper variety of groups  
${\cal G}$ that  the unitriangular automorphism group  $U_n$ is linear  if and only if the  relatively free group $G_{n-1}$  is virtually nilpotent. More exactly  (for $n \geq 3$): if $G_{n-1}$ is virtually nilpotent, then $U_n$ admits a faithful matrix representation over integers  ${\mathbb{Z}}.$ It was also shown in \cite{Erofeev} that if  $n \geq 3$ and  $G_{n-1}$ is nilpotent then    $U_n$ is nilpotent too.

Now let $C_n$ be an arbitrary relatively free algebra of rank  $n$ with  set of free generating elements  $\{x_1, ... , x_n\}.$ For    $m \leq n$    $C_m$ can be considered as  subalgebra of        $C_n$ generated by  $x_1, ... , x_m.$ An automorphism   $\psi $ of $C_n$ is called {\it unitriangular}
w.r.t. the given set of free generating elements if it is defined by map of the form:

\begin{equation}
\label{eq:2}
\psi : x_1 \mapsto x_1, x_i \mapsto x_i + u_i \   {\textrm{для}} \ i = 2, ... , n,
\end{equation}  

\noindent
where $u_i = u_i(x_1, ... , x_{i-1})$ belongs to  $C_{i-1}.$ Let $U_n$ denote a subgroup of the automorphism group Aut$(C_n)$ of $C_n,$ consisting of  all unitriangular  automorphisms. As abstract group  $U_n$  does not depend from a chosen set of free generating elements of $C_n.$  

The author, Chirkov and Shevelin    \cite{Romankov} proved that, for a free Lie (free associative, absolutely free, polynomial) algebra $C_n$ of rank $n \geq 4$ over a field of zero characteristic, the unitriangular automorphism group  $U_n$  is not linear.  Then  the following papers \cite{Sosnovskii}, \cite{Kabanov} presented descriptions of the hypercentral series of groups  $U_n$ corresponding to polynomial and free metabelian Lie algebras, respectively. By these results  $U_n$ are not linear for $n \geq 3.$ By \cite{Bardakov}, for $n \geq 3,$ the unitriangular automorphism group $U_n$ is not linear in case of polynomial algebra and in case of free associative algebra.  
 By  \cite{Romankov1}  for each relatively free algebra $C_n$ the group $U_n$ is locally nilpotent, thus it is linear.

\section{The method of Formanek and Procesi.}

\medskip
\noindent
Let   $G$  be any group, and let ${\cal H}(G)$ denote the following HNN-extension of $G \times G$:

\begin{equation}
\label{eq:3}
{\cal H}(G) = <G \times G, t : t(g,g)t^{-1} = (1,g), g \in G>.
\end{equation}

\noindent
\begin{theorem}
{ (Formanek, Procesi \cite{Formanek}).} \textit{Let $\rho $ be a linear representation of ${\cal H}(G).$ Then the image of  $G \times \{1\}$ has a subgroup of  finite index with nilpotent derived subgroup, i.e,  is nilpotent-by-abelian-by-finite.} 
\end{theorem}

\medskip
\noindent
\begin{theorem}
{(Brendle, Hamidi-Tehrani \cite{Brendle}).}
 \textit{Let     $N$   be a normal subgroup of ${\cal H}(G)$   such that the image of    
$G \times \{1\}$ in  ${\cal H}(G)/N$   is not nilpotent-by-abelian-by-finite. Then  ${\cal H}(G)/N$  is not linear.}
\end{theorem}

In \cite{Brendle} a group of the type described in Theorem 3.2 is called  a {\it Formanek and Procesi group}, or {\it FP-group} for short.

\section{Proof of  Theorem 1.1.}

\medskip
\noindent
{\bf 3.1.} For   $n \leq 3,$       $U_n$  is linear.

\medskip
\noindent
{\it Proof:}  
Since $U_1$  is trivial and $U_2$ is infinite cyclic the statement is obvious for $n = 1, 2.$

Let $n = 3.$ By\cite{Erofeev}  $U_3$ is generated by automorphisms $\lambda_{2,1}, \lambda_{3,1}, \lambda_{3,2}.$
Recall that $\lambda_{i,j}$ maps $f_i$ to $f_jf_i,$ and fixes all other basic elements. This is
applicable for any group $U_n.$ The automorphisms $\lambda_{3,1}$   
$\lambda_{3,2}$ generate in $U_3$ a normal free subgroup $F_2.$ The automorphism  $\lambda_{2,1}$ acts as follows:

\begin{equation}
\label{eq:4}
\lambda_{2,1}^{-1}\lambda_{3,1}\lambda_{2,1} = \lambda_{3,1}, \ 
\lambda_{2,1}^{-1}\lambda_{3,2}\lambda_{2,1} = \lambda_{3,1}\lambda_{3,2}.
\end{equation} 

Now we'll show that    $U_3$  is isomorphic to a subgroup of  Aut$(F_2).$ Let 
$\tau_1$   and  $\tau_2$ denote inner automorphisms of $F_2$ corresponding to   $f_1$   and  $f_2$  respectively. This means that any element  
$g$ of  $F_2$ maps by  $\tau_i (i=1,2)$ to
 $f_i^{-1}gf_i.$  Let  $\sigma_{2,1}\in $ Aut$(F_2)$ fixes     $f_1$  and maps  $f_2$ to 
$f_1f_2.$ Obviously,   $F_2=$ gp($\tau_1, \tau_2)$ is a free group of rank  $2.$  It is a normal subgroup of   $V_3 =$ gp$(\tau_1, \tau_2, \sigma_{2,1}).$ A quotient $V_3/F_2$ is the infinite cyclic generated by the image of  $\sigma_{2,1}.$ The corresponding action is determined by: 

\begin{equation}
\label{eq:5}
\sigma_{2,1}^{-1}\tau_{1}\sigma_{2,1} = \tau_{1}, \ 
\sigma_{2,1}^{-1}\tau_{2}\sigma_{2,1} = \tau_{1}\tau_{2}.
\end{equation} 

Thus,   $U_3$ and  $V_3$ are both infinite cyclic extensions of  $F_2.$ By
 (4) and (5) we conclude that $\alpha : U_3 \rightarrow V_3$  defined as:

\begin{equation}
\label{eq:6}
\alpha  : \lambda_{3, j} \mapsto \tau_{j} \ {\textrm{для}} \  j = 1, 2, \  
\lambda_{2,1} \mapsto \sigma_{2,1},
\end{equation} 

\noindent
is isomorphism. Since  $V_3$ is a subgroup of   Aut$(F_2),$  which is linear by \cite{Dyer} and \cite{Krammer},      $U_3$ is also linear. 

\hfill $\Box$

\medskip
\noindent
{\bf 3.2.}  For  $n \geq 4,$     $U_n$  is not linear. 

\medskip
\noindent
{\it Proof:} For $n \geq m,$    $U_n$ has a subgroup that is isomorphic to  $U_m.$ Elements of this subgroup act naturally to $f_1, ... , f_m$ and fix elements $f_{m+1}, ... , f_n.$  So, we just have to prove that  $U_4$  is not linear. 

By Theorem 3.2 it will be enough to find a subgroup $H$ of  $U_4$ that is isomorphic to a quotient    ${\cal H}(F_2)/N,$  where ${\cal H}(F_2)$ is given by  (3), such that the image of  $G \times \{1\}$  in  ${\cal H}(G)/N$   is not nilpotent-by-abelian-by-finite.   

 There are two commuting elementwise subgroups  of $U_4$ each of them is isomorphic to  $F_2.$ Namely, there are  gp$(\lambda_{3,1}, \lambda_{3,2})$ and gp$(\lambda_{4,1}, \lambda_{4,2}).$  Consider them as two copies of $F_2$ via isomorphism defined by map $\lambda_{3,1} \mapsto \lambda_{4,1}, \lambda_{3, 2} \mapsto \lambda_{4,2}.$  Thus we have a subgroup $F_2 \times F_2$  of  $U_4.$  Easily to check that:

\begin{equation}
\label{eq:7}
\lambda_{4,3}^{-1}\lambda_{3,j}\lambda_{4,j}\lambda_{4,3} = \lambda_{3,j}, \  j = 1, 2.
\end{equation}

 By (3) and (7) we conclude that  $H$ is a homomorphic image of  ${\cal H}(F_2)$ such that the subgroup   $F_2 \times F_2$ of ${\cal H}(F_2)$ maps isomorphically to the just constructed subgroup of the same type of $U_4.$  
The image of   $t$ is $\lambda_{4,3}.$ Hence, $H \simeq {\cal H}(F_2)/N,$  where   $N$ is the kernel of this homomorphism. By Theorem 3.2  $H,$ and so $U_4,$ is not linear. 

\hfill $\Box$

\medskip
\noindent\textbf{Remark  1.} \emph{In fact we proved that subgroup
$W_4=$ gp$(\lambda_{3,j}, \lambda_{4,l} \  : j = 1, 2; l = 1, 2, 3)$  of  $U_4$ is not linear. In \cite{Erofeev} we noted that  $[\lambda_{i,j}, \lambda_{j,k}] = \lambda_{i,k}.$ Here  commutator  $[g, f]$ means    $gfg^{-1}f^{-1}.$ Any group  $U_n$ is generated by the elements  $\lambda_{i,j},$ for $j < i \leq n$  (see \cite{Erofeev}).  It follows that $W_4$  is the derived subgroup       $U_4'$   of $U_4.$ Hence, we proved that the derived subgroup (the second member $\gamma_2U_4$ of the low central series)   of $U_4$
is not linear.   This subgroup  $\gamma_2U_4$  has also characterized in $U_4$ as the stabilizer of
 $f_1.$ In general case, for  $n \geq 4,$  member  $\gamma_{n-2}U_n$ coincides with the elementwise stabilizer of  $\{f_1, ..., f_{n-3}\}.$ 
Easily to see that the derived subgroup  $U_4'$ can be embedded into  $\gamma_{n-2}U_n.$  Hence, for every 
$n \geq 4,$ a member  $\gamma_{n-2}U_n$ of the low central series of $U_n$ is not linear. This statement is more strong than the statement of Theorem 1.1 about nonlinearity of $U_n$ for $n \geq 4.$}    

\medskip
\noindent\textbf{Remark 2.} \emph{In  \cite{Bardakov1} an explicit faithful representation of   Aut$(F_2)$  in  GL$_{12}({\mathbb{Z}}[t^{\pm 1}, q^{\pm 1}])$  is given. Here 
${\mathbb{Z}}[t^{\pm 1}, q^{\pm 1}]$ is a  Laurent polynomial ring. Hence, $U_3$ has a faithful matrix representation  over ${\mathbb{Z}}[t^{\pm 1}, q^{\pm 1}].$ }

 \emph{Also note that  $U_3$ can be presented by $<a, b : [[a, b], b]] =1>,$ where  $a$ corresponds to   $\lambda_{3,2},$ and  $b$ corresponds to  $\lambda_{2,1}.$  By terminology of  \cite{Baumslag}  $U_3$ is the first non commutative member in a series of Hydra groups  $H_k = <a, b : [ ...  [ [a, b], b], ... , b] = 1>, k \geq 1,$ where  commutator has $k$ entries of $b.$ In general, Hydra groups were introduced in \cite{Dison}.
It was shown in \cite{Baumslag} that Hydra groups of such  form are residually torsion-free nilpotent. It seems interesting to study  their linearity.} 
 
\medskip
\noindent\textbf{Remark 3.} \emph{By   \cite{Bajorska}  a group $G$ is called {\it locally graded} if every nontrivial finitely generated subgroup of $G$ has a proper subgroup of finite index. This class contains, for example, all locally solvable and all residually finite groups. Let  $G_2$ be the relatively free group of rank $2$ in the variety var$(G)$ generated by  $G.$ Suppose that the derived subgroup  $G_2'$ is finitely generated. Then by  \cite{Bajorska}  $G$ is virtually nilpotent.  }  

\emph{ Let ${\cal G}$ be a variety consisting of locally graded groups. Obviously ${\cal G}$ is a proper variety.  Suppose that  $U_3$ is linear. Then every group  $U_n$ is linear. Indeed,
 by \cite{Erofeev}  $G_2$ is virtually nilpotent. It follows that $G_2'$ is finitely generated. Any group $G_{n-1}$ generates a subvariety  var$(G_{n-1})$  of ${\cal G}.$ The relatively free group of rank $2$ in this subvariety is a homomorphic image of $G_2,$ 
and so has finitely generated derived subgroup. Then by  
\cite{Bajorska}  $G_{n-1}$ is virtually nilpotent. It follows by \cite{Erofeev} that  $U_n$ is linear. Thus, the linearity of   $U_3$  implies the linearity of  $U_n$ for every  $n \geq 4.$ Moreover,  ${\cal G}$ should be virtually nilpotent.}

\emph{ We see by Theorem 1.1 that just presented statement, that the linearity of $U_3$ implies the linearity of $U_n$ for all $n \geq 4,$ is not true for the variety of all groups.  Likely, it is also non-true  for some proper varieties of groups. As  candidates to such varieties we can consider the varieties of groups generated by the famous Golod groups. We conjecture that for every  $m \geq 3$ there is a variety  ${\cal G}_m$ such that  
the groups $U_n$ are linear if and only if   $n \leq m.$ }

\end{document}